\documentclass[twoside]{article}
\begin{document}
\mathchardef\sym="3218
\markboth{IMRE MAJOR}{STRATIFICATION BY TRAJECTORY SPACES}
\author{Imre Major}
\title{On the stratification of a compact 3-manifold by the trajectroy spaces of a Morse-Smale flow}
\date{}
\maketitle

$\hspace{2.4cm}$Dedicated to Professor Paul Nevai$\vspace{1.3cm}$

A B S T R A C T $\vspace{.3cm}$

We consider a Morse function $f$ and a Morse-Smale gradient-like vector field $X$ on a compact connected oriented 3-manifold
$M$ such that $f$ has only one critical point of index $3$. Based on Laudenbach's ideas \cite{L}, we will show that the flow of $X$
can be isotoped into one so that the trajectory spaces of the new flow provide a stratification for $M$. We will construct "natural"
tubular neighborhoods about each given trajectory space of the new flow such that these neighborhoods are stratified by open subsets of
trajectory spaces that co-bound the given one. In connection with this we introduce the concept of {\it conic stratification} of a
manifold and point out that this is the appropriate condition the stratification of $M$ by trajectory spaces should be required to satisfy.


\section{Introduction, statement of results.}

Let $M$ be a compact oriented 3-manifold (possibly with boundary) equipped with an $f:M\to{\bf R}$ self-indexing Morse function and let $X$
be a Morse-Smale gradientlike vector field for $f$ (we presume $\partial M=f^{-1}(0)$ whenever it is non-empty). We call $(f,X)$ a Morse-pair.
(Basics of Morse Theory can be found in \cite{Mi1} or \cite{Mi}.)

By Morse's Cancellation Theorem (see \cite{Mi} Theorem 5.4) we can suppose that there is only one critical point of index $3$ (denoted by
$o$). Let$$\{x_0,...,x_K\},\hspace{1cm}(\{y_1,...,y_L\})$$be the sets of critical points of index two (one) respectively, thus the {\it
critical set} of  $f$ is $${\rm Cr}(f)=\{o\ ,x_1,...,x_K,\ y_1,...,y_L,\ b_1,...,b_J\}$$where $ b_1,...,b_J$ are the critical points of index $0$ (with
$J=0$ when $\partial M\not=\emptyset).\vspace{.3cm}$

We fix a Morse-chart $(U_o,\eta_o)$ about the top critical point $o$ and Morse charts$$\eta_{x_k}:U_{x_k}\to{\bf R}^3\hspace{.3cm}
(k=1,...,K),\hspace{.3cm}\eta_{y_l}:U_{y_l}\to{\bf R}^3\hspace{.3cm}(l=1,...,L)$$about the rest of  the critical points (except for the $b_j$'-s
at which we don't need such chart). We suppose that all the  images$$\eta_o(U_o)=\eta_{x_k}(U_{x_k})=\eta_{y_l}(U_{y_l})$$are 3-balls about the
origin with radius $2\delta$. Let \begin{eqnarray*}S^\downarrow_k&:=&\eta_{x_k}^{-1}(S_{\delta}\cap({\bf R}^2\times\{0\}))\\ S^\uparrow_l&:=
&\eta_{y_l}^{-1}(S_{\delta}\cap(\{0\}\times{\bf R}^2))\end{eqnarray*}be the pre-images of the $\delta$-circles in the negative (positive)
subspace of the Hessian of  $f\circ\eta_{x_k}^{-1}\hspace{.2cm}(f\circ\eta_{y_l}^{-1})$ at the origin of ${\bf R}^3$. (We used the canonical
identification $T_{\bf 0}{\bf R}^3\equiv{\bf R}^3$ for the tangent space at the origin.) We call $S_k^\downarrow$ the {\it outbound  circle }
(of the flow) at critical point $x_k$ ($S_l^\uparrow$ is the {\it inbound circle} at $y_l$, respectively.)

Based on the identifications $U_{x_k}\mathop{\equiv}\limits^{\eta_{x_k}}B_{2\delta}$ etc., we will treat the Morse-domains $U_o$,
$U_{x_k}$, $U_{y_l}$ as if they were 3-balls (in other words, we will use the Euclidean structure on each Morse domain without
any further reference to the Morse coordinate systems).$\vspace{.3cm}$

{\noindent\it Notations:}  $\lambda_p$ stands for  the trajectory of vector field $X$ through point $p\in M$. Observe that each trajectory
$\lambda$ of $X$ is defined on the entire real line and $$\mathop{\lim}\limits_{t\to-\infty}\lambda(t)\hspace{1cm}\mathop{\lim}\limits_{t
\to+\infty}\lambda(t)$$ are critical points. They are called the {\it initial} and the {\it terminal} point of $\lambda$, respectively. We
say that a trajectory $\lambda$ {\it connects} its initial and terminal points. The {\it index difference} of $\lambda$ is the difference of
indices at its initial and terminal points. Let $W_l^\uparrow\ (W_k^\downarrow)$ denote the stable (unstable) subsets of vector field $X$
at critical point $y_l\ (x_k)$ respectively. Similarly, $W_o^\downarrow$ denotes the unstable submanifold of  $X$ at critical point $o$. Then
we can define {\it trajectory spaces}$$W_{ol}:=W_o^\downarrow\cap W_l^\uparrow,\hspace{1cm}W_{kl}:=W_k^\downarrow\cap W_l^\uparrow$$etc.
Finally, for an arbitrary subset $N\subset M$ let $\cal N$ denote the set of trajectories ${\cal N}=\{\lambda_p\ |\ p\in{\cal N}\}$ (e.g.
${\cal W}_k^\downarrow$ is the set of trajectories that emanate from critical pt. $x_k,\ {\cal W}_{kl}$ is the set of trajectories that connect
critical points $x_k$ and $y_l$, etc.) Fix an orientation on each unstable set $W_k^\downarrow$ at critical point $x_k\ (k=1,\dots,K)$ and
stable set $W_l^\uparrow$ at critical point $y_l\ (l=1,\dots,L)$. Observe that trajectory space ${\cal W}_{ok}$ consists of two trajectories
$\mu_k^-$ and $\mu_k^+$ where a pair of positively oriented vectors of $T_{x_k}W_k^\downarrow$ together with the orientation of $\mu_k^+$ are
positively oriented. Similarly, the unstable manifold of vector field $X$ at critical point $y_l$ is the union of trajectories $\nu_l^-\cup
\nu_l^+$.

The orientation on $W_l^\uparrow$ provides a cyclic ordering on the set of trajectories that initiate at a critical point of index $2$ and
terminate at fixed critical point $y_l$. We list this set as$$\lambda_l^1,\dots,\lambda_l^{s_l}$$(indexing is understood in the cyclic sense,
i.e.  modulo $s_l$).

It is known that when $\partial M=\emptyset$ then a CW-decomposition of $M$ can be associated to a Morse-pair $(f,X)$ with sole 3-cell $W_o
^\downarrow\cup\{o\}$, 2-cells$$W_k^\downarrow\cup\{x_k\}\hspace{.5cm}(k=1,\dots,K),$$1-cells$$\nu_l^-\cup\nu_l^+\cup\{y_l\}\hspace{.5cm}(l=1,
\dots,L)$$and $0$ cells: $\{b_1,\dots,b_J\}$. This known result can easily be re-proven by the techniques we provide.\vspace{.3cm}

Now we will describe the conic stratification structure induced by the Morse-pair. The concept of conic stratification works in all finite
dimensions. One of our future plans is to generalize the theorem given below to arbitrary finite dimensions and also, to the $G$-case (i.e.
to the case when a compact group $G$ acts on $M$, the Morse function $f$ is invariant and $X$ is equivariant with respect to this action).
(Pre)stratifications are discussed in detail in Mather \cite {Ma} and Verona \cite{V1}, \cite{V2}.\vspace{.3cm}

{\noindent\bf Definition 1.} A {\it prestratification} of a differentiable manifold $M$ is a locally finite partition $\{M_\alpha\}_{\alpha
\in{\cal A}}$ of $M$ into submanifolds $M_\alpha$ (called {\it strata} and indexed by $\alpha\in\cal A$) such that the following {\it frontier
condition} holds:\vspace{.3cm}

The {\it frontier} Fr($M_\alpha)=\overline M_\alpha\setminus M_\alpha$ of  $M_\alpha\ (\alpha\in{\cal A})$ is a union of

certain strata.$\vspace{.3cm}$

{\noindent\bf Example:} Given a Morse-pair $(f,X)$, the critical points of $f$ and the trajectory spaces provide a prestratification of $M.\vspace{.3cm}$

We use notation $\alpha'\prec\alpha$ whenever $M_{\alpha'}\subseteq $Fr$(M_\alpha)$. Then the relation "$\prec$" is a partial order
on the index set $\cal A$. Note that for a given prestratification $\{M_\alpha\}_{\alpha\in{\cal A}}$ of $M$ and an open subset $U\subset M$ partition
$\{M_\alpha\cap U\}_{\alpha\in{\cal A}}$ is a prestratification of $U$ called the {\it induced prestratification}. In order to simplify
notation, the indexing will always be understood that only the non-empty constituents are taken into account (i.e. subsets in a partition,
trajectory spaces between critical points, etc), repeated indices only once.\vspace{.3cm}

Out of the available ways, we choose to define a {\it tubular neighborhood} of an oriented submanifold $N\subset M$ as a sixtuple $(E,\rho,N,
O,\epsilon,\xi)$ where:$\vspace{.2cm}$

$E\mathop{\longrightarrow}\limits^\rho N$ is a smooth oriented vector bundle with an Hermitian

structure $\langle,\rangle.\vspace{.2cm}$

$O\subset M$ is an open neighborhood of $N$, called the {\it tube} about $N.\vspace{.2cm}$

$\epsilon:\overline N\to[0,\infty)$ is a continuous function, $\epsilon^{-1}(0)=\overline N\setminus N.\vspace{.2cm}$

$\xi:E(\epsilon)\to O\hspace{.3cm}\ (\xi|_N=$id$_N)$ is a diffeomorphism from  the $\epsilon$-disc

bundle of $E$ onto $O$ so that its differential preserves the orientation

along the zero section $N\subset E.\vspace{.2cm}$

{\noindent\bf Definition 2.} A {\it conic prestratification of a vector bundle} $E\mathop{\longrightarrow}\limits^\pi B$ consists of:$\vspace{
.2cm}$

(i) A filtration $B=E^0\subseteq E^1\subseteq...\subseteq E^n=E$ so that the $E^i$-s are

invariant under multiplication by {\underline{positive}} real numbers and the

partition $\{B\}\cup\{E^i\setminus E^{i-1}|i=1,...,n\}$ is a prestratification of $E.\vspace{.2cm}$

(ii) For each $b\in B$ there is an open neighborhood $U\subset B$ about $b$

and a local trivialization $t_b:\pi^{-1}(U)\to U\times\pi^{-1}(b)$ which is com-

patible with the filtration (i.e. for which the following holds:$$t_b(\pi^{-1}(U)\cap E^i)=U\times(\pi^{-1}(b)\cap E^i)$$

{\noindent\bf Definition 3.} Let $E\mathop{\longrightarrow}\limits^\pi B$ be an Hermitian vector bundle. Suppose that $B\subset M$ is a
submanifold of the compact Riemannian manifold $(M,{\bf g})$ and let

{\noindent}$\epsilon:\overline{B}\to[0,\infty)$ be a continuous function for which $\epsilon^{-1}(0)=$Bd$(B)$ holds. We say that the $\epsilon
$-disc bundle $E(\epsilon)$ is {\it thick} if $\exists\ a>0$ such that$$\epsilon(x)>ad_{\bf g}(x,{\rm Bd}(B))$$where $d_{\bf g}$ is the
Riemann-distance on $M$ induced by metric $\bf g$.\vspace{.3cm}

{\noindent\it Remark:} For a thick disc bundle one can always find positive reals $a'$ and $\delta'$ such that $\epsilon(x)>a'\delta'$ for
each $x\in B$ with $d_{\bf g}(x,{\rm Bd}(B))>\delta'$. By this reason, in the sequel we will presume that the continuous function $\epsilon$
associated to a thick  neighbourhood is of the form \begin{eqnarray*}\epsilon(x)&=&ad_{\bf g}(x,{\rm Bd}(B))\hspace{.3cm}{\rm for}\
d_{\bf  g}(x,{\rm Bd}(B))\leq\delta)\\ \epsilon(x)&=&a\delta\hspace{2.45cm}{\rm for}\ d_{\bf g}(x,{\rm Bd}(B))>\delta\end {eqnarray*}

Fix a Reimannian metric $\bf g$ on $M$ which restricts to the Eucledian metric on the Morse-domains. A tubular neighbourhood of a submanifold
$N\subset M$ is {\it thick} if the function $\epsilon$ defines a thick disc bundle. It is called {\it  tangential} if $E\subset T_NM$ is a
subbundle and$$v={d\over dt}|_{t=0}\xi(tv)\hspace{.3cm}(v\in E)$$holds. In the sequel we will consider tangential tubular neighborhoods only. A tangential
tubular neighborhood is {\it geodesic} if $\xi$ is the restriction of the exponential map (of the Live-Civita connection of metric $\bf g$)
to subbundle $E\subset T_NM$. For a tubular neighborhood $(E,\rho,N,O,\epsilon,\xi)$ the map$$r:O\to N,\hspace{.6cm} r:=\rho\circ\xi^{-1}$$is
called the {\it associated retraction} while$$d_N:O\to[0,\infty),\hspace{.6cm}d_N(p):=\sqrt{{\bf g}(\xi^{-1}(p),\xi^{-1}(p))}$$is the {\it
distance from submanifold} $N.\vspace{.3cm}$

The definition of a "stratified tubular neighborhood" can now be provided as follows:\vspace{.2cm}

{\noindent\bf Definition 4.} The prestratification $\{M_\alpha\}_{\alpha\in{\cal A}}$ of Riemannian manifold $M$ is called {\it conic
prestratification} if about each strata $M_\alpha$ there is given a tubular neighbourhood {\noindent}$(E_\alpha,\rho_\alpha,M_\alpha,O_\alpha,
\epsilon_\alpha,\xi_\alpha)$ such that the vector bundle $(E_\alpha,\rho_\alpha,M_\alpha)$ is conically prestratified by subbundles$$\{E
_\alpha^\beta\ |\ \beta\in{\cal A},\ \alpha\prec\beta\},$$the disc bundle $E_\alpha(\epsilon_\alpha)$ is thick and diffeomorphism $\xi_\alpha
:E_\alpha(\epsilon_\alpha)\to O_\alpha$ preserves strata with respect to the induced prestratification on $O_\alpha$ (i.e.$$\xi_\alpha(E
_\alpha^\beta(\epsilon_\alpha))=O_\alpha\cap M_\beta)$$

{\noindent\it Remark:} Local finiteness of the family $\{M_\alpha\}_{\alpha\in{\cal A}}$ implies that it is possible to choose the open
neighborhoods $O_\alpha$ so that:

(i) $O_\alpha\cap M_\beta\not=\emptyset$ iff $\alpha\prec\beta$.

(ii) $O_\alpha\cap O_\beta\not=\emptyset$ iff $\alpha$ and $\beta$ can be compared w.r.t. $\prec.\vspace{.3cm}$

{\noindent\bf Definition 5.} With notations from the previous definition, a conic prestratification is called {\it conic stratification} if
compatibility conditions\begin{eqnarray}r_\beta(p)\in O_\alpha,\hspace{1cm}r_\alpha\circ r_\beta(p)=r_\alpha(p)\hspace{1cm}(p\in O_\alpha
\cap  O_\beta,\ \alpha\prec\beta)\end{eqnarray} hold for the associated retractions.$\vspace{.3cm}$

{\noindent\it Remark:} Conditions in (1) are the two Thom-Mather control conditions combined together, the first one in a somewhat weakened
form. A {\it stratification} (in the usual sense) is gotten from a prestratification by stipulating that for each strata a retraction and a
height function is associated so that (1) holds moreover the height function remains invariant under the retractions assigned to higher
dimensional strata (See Verona \cite{V2} pg. 3 or Mather \cite{Ma} pg. 33). This general concept does not involve the bundle property.
Laudenbach (\cite{L} Proposition 2.) proves that (in any finite dimension) the trajectory spaces of a Morse-Smale flow consitute (what he
calls) a {\it submanifold with comic singularities} (abbreviated smcs). He gives an inductive definition of smcs similar to Definition 4.
(without requiring formula (1)). In this article we will describe a way to isotop the flow so that a conic stratification can be constructed
for the trajectory spaces of the new flow (in a natural way). An application of the inverse isotopy then provides conic stratification structure
for trajectory spaces of the original flow. Thus we will give an alternative discussion of Laudenbach's Proposition 2., ensuring that
formula (1) holds this time. It is also important to have a way to gauge the size of the domains (called {\it tubes about strata}) on which
formula (1) holds (i.e. the thickness, in our terminology). We will return to the major difficulties of generalizing our treatment to arbitrary finite dimesion
at the end of this paper.\vspace{.3cm}

Let $S_o:=\eta_o^{-1}(S_{\delta})\subset U_o$ denote the outbound sphere at the top critical point $o$ and let $S_{ol}=S_o\cap W^\uparrow
_l\hspace{.3cm}(l=1,\dots L)$. Then partition\begin{eqnarray}\{S_o\setminus[\mathop{\cup}\limits_{k=1}^KW_k^\uparrow\cup\mathop{\cup}\limits
_{l=1}^LS_{ol}],\ S_o\cap W_k^\uparrow,\ S_{ol}\ |\ k=1,\dots K,\ l=1,\dots,L\}\end{eqnarray}is a prestratification of the outbound sphere
$S_o$.\vspace{.4cm}

{\noindent\bf Theorem A:} Given a Morse-pair $(f,X)$, there exists an isotopy $H:M\times I\to M$, stationary in a neighborhood of critical
points such that:$\vspace{.2cm}$

(i) $H_0={\rm id}_M$.\vspace{.2cm}

(ii) For the new Morse-pair $(f',X')$ with $f'=f\circ H_1,\hspace{.3cm}X'=$d$H_1(X)$ the

geodesic tubular neigborhoods of the $0$ and $1$-dimensional starta in (2) provide

a conic stratification of sphere $S_0$.\vspace{.2cm}

(iii) The above conic stratificaion of $S_0$ extends naturally (via the flow) into a

conic stratification of $M$ by the trajectory spaces of the new flow.$\vspace{.5cm}$

{\noindent\bf Theorem B:} There exist:\vspace{.3cm}

(i) $L$ disjoint differentiable curves on the unit 2-sphere $S=$Bd($B^3$) such that

 the connected components are circles or open curves.\vspace{.2cm}

(ii) A surjection $\psi:\overline{B^3}\to M\ \ni\ \psi|_{B^3}:B^3\to W_o^\downarrow$ is a diffeomorphism that

takes the open cone over one of the above curves (with the origin as vertex)

diffeomorphically onto one of the trajectory spaces $W_{ol}$.\vspace{.2cm}

(iii) Surjection $\psi$ serves as a characteristic map for 3-cell $W_o^\downarrow$.\vspace{.4cm}

The conic nature of an arbitrary prestratification is discovered by Verona (\cite{V1}, Propositions 2.1. and 2.6.) but in this general setting
we loose local triviality. We will discuss how his ideas fit into the picture in the end of this paper, after we will have developed the
necessary tools to provide such discussion.\vspace{.4cm}

{\section{Proof of  the theorems}

{\noindent}For a pair of points $p,q\in{\bf  R}^3$ the (open) straight section between $p$ and $q$ is denoted by $\overline{pq}.\vspace{.2cm}$

{\noindent\bf Definition 6.} ({\it Standard 3-cell}): Let $B\subset{\bf  R}^3$ be the open  ball about the origin with radius $3$ and
define the map\begin{eqnarray*}\phi_o:B&\hookrightarrow&M,\hspace{.5cm} \phi_o(v)=\lambda_{p_v}\cap f^{-1}(3-\|v\|)\hspace{.5cm}(\|v\|\not=0)
\\ \phi_o({\bf 0})&=&o\end{eqnarray*}where $p_v:=\eta_o^{-1}({\delta v\over3})\vspace{.3cm}$

{\noindent\it Observe:} that the domain $B_o $ of  map $\phi_o$ is a union of open (straight) sections that initiate at the origin  and the
origin itself (thus $B_o$ is a star-shaped open subset of  $B$) moreover $\phi_o$  maps each such open section onto a trajectory of the vector
field $X$ such that$$\phi_o:B_o\to W_o^\downarrow\cup\{o\}$$ is a diffeomorphism. The open sections in $B_o$ (that emanate from the origin)
can be classified as:$\vspace{.4cm}$

{\bf I.}  sections of length 1: these are mapped onto trajectories of  $X$ that
terminate

at a critical point of index 2 (i.e. onto one of  $\mu_k^\pm\hspace{.2cm}(k=1,...,K)).\vspace{.3cm}$

{\bf  II.} sections of length 2: these are mapped to trajectories that terminate

at a critical point of index one. $\vspace{.3cm}$

{\bf  III.}  the length of the rest of the sections is 3 and they are mapped onto

trajectories that terminate at one of critical points $b_1,..,b_J$ (or at $\partial M).\vspace{.3cm}$

{\noindent}The set of  trajectories ${\cal W}_o^\downarrow$ is in a bijective correspondence with sphere $S:=$Bd$(B)$ (rays that emanate from
the origin $\bf  0$ in ${\bf R}^3$ are in bijective correspondence with points of  $S$ and each such ray contains a unique section that is
taken to a trajectory of $X$ by map  $\phi_o$). We call these sections (and the rays that contain them) {\it sections of type} {\bf I., II.}
and {\bf III.}, according to the above classification and say that points of $S$ {\it represent} trajectories in ${\cal W}_o^\downarrow.
\vspace{.2cm}$

{\noindent}Sections of type {\bf  I.} come in pairs that correspond to the pairs of trajectories $(\mu_k^-,\mu_k^+)$ thus the pairs of rays
that contain these sections meet the sphere $S$ in a pair of points $(x_k^-,x_k^+)$. We say, that points $(x_k^-,x_k^+)$ are {\it opposite}
to each other.

$\vspace{.2cm}$

{\noindent}For type {\bf II.} sections consider the trajectory space $W_{ol}$ between critical points $o$ and $y_l$. This trajectory space
can be gotten from the ascending manifold $W_l^\uparrow$ by deleting all trajectories that initiate at a critical point of index 2, so the
connected components of  $W_{ol}$ are represented by oriented curves $\alpha_l^1,...,\alpha_l^{s_l}\subset S$. Each point of $\alpha_l^i$
represents a trajectory in ${\cal W}_{ol}$ which lies between trajectories $$\lambda_l^i,\hspace{.2cm}\lambda_l^{i+1}\in{\cal W}_l^\uparrow
\hspace{1cm}(i=1,...,s_l)$$The union of  trajectories represented by curve $\alpha_l^i$ is bounded by two joinig pairs of trajectories$$
(\mu^\varepsilon_k,\lambda^i_l),\hspace{.3cm}(\mu^{\varepsilon'}_{k'},\lambda^{i+1}_l))$$where $\varepsilon,\ \varepsilon'\in\{-,+\}$ so curve
$\alpha_l^i$ connects points $x_k^\varepsilon$ and $x_{k'}^{\varepsilon'}\ (k=k'$ is allowed). Observe that the predecessor, $\alpha_l^{i-1}$
terminates at the {\underline{opposite}} of the initial point of  $\alpha_l^i$. To see this, note that trajectory $\lambda_l^i$ bounds two
connected components of the trajectory space $W_{ol}$ (represented by curves $\alpha_l^{i-1}$ and $\alpha_l^i$, respectively). The intersection
of these two components with the Morse domain  $U_{x_k}$ fall into different components of  $U_{x_k}\setminus(W_k^\downarrow\cup\{x_k\})$,
consequently one of the components is bounded by trajectory $\mu_k^+$ and the other one is bounded by $\mu_k^-$.

{\noindent}It follows,  that  the closure of curve $$\alpha_l:=\mathop{\cup}\limits_{i=1}^{s_l}\alpha_l^i$$is {\underline{not connected}}.$\vspace{.3cm}$

{\noindent Observe} that map $\phi_o:B_o\to W_o^\downarrow$ is a diffeomorphism, thus $W_o^\downarrow$ is a 3-cell and $$M\setminus\phi_o(B_o)
={\big(}\mathop{\cup}\limits_{k=1}^K W_k^\downarrow\cup\{x_k\}\big)\cup\mathop{\cup}\limits_{l=1}^L \overline{\nu_l^-\cup\nu_l^+}$$ while
$B\setminus B_o$ is a union of (left-closed) sections of rays through points of set $\mathop{\cup}\limits_{l=1}^L\overline{\alpha_l}$. We call
$B\setminus B_o$ the {\it singular part} of ball $B$. It's clear that in its present form  $\phi_o$ can not be extended to the boundary of
ball $B_o$ (this boundary is not even a 2-sphere). In what follows we will modify $\phi_o$ into a characteristic map for a 3-cell. In order to
attain this, we will construct a   homotopy \begin{eqnarray}H:B_o\times[0,2]\to B_o\end{eqnarray}such that:$\vspace{.2cm}$

(i) $H_0=$id$_{B_o}.\vspace{.3cm}$

(ii) Each $H_t$ is a diffeomorphism onto an open subset of  $B_o\hspace{.2cm}(0\leq t\leq 2).\vspace{.3cm}$

(iii)  The boundary  $S^*:={\rm Bd}(H_2(B_o))$ (taken in $\overline B$) is homeomorphic to

the 2-sphere. $\vspace{.3cm}$

{\noindent}Roughly speaking, we will "push in" the star-shaped set $B_o$ along  its singular part $B\setminus B_o$ in two steps: first, along
the singular part of rays of type {\bf  I.} and then along the rest of the singular set $B\setminus B_o$. To perform the above procedure, first
we need to alter the flow so that it becomes "nice" (we call such an altered, nice flow an {\it aligned} flow).\vspace{.4cm}

{\noindent\underline{Alignment of the flow:}} we will alter the flow in the neighborhood of each trajectory with index difference one (finitely
many trajectories altogether). As the process is the same for all such trajectories, we will describe it for a trajectory  $\lambda\in{\cal
W}_{kl}$ for some fixed pair of indices $(k,l)$. Let \begin{eqnarray*}S_k&:=&\eta_{x_k}^{-1}\ (S_{\delta})\\  S_l&:=&\eta_{y_l}^{-1}(S_{\delta}
)\\ p&:=&\lambda\cap S_k^\downarrow\\  q&:=&\lambda\cap S_l^\uparrow\\  P_p&:=&T_pS_k\\  P_q&:=&T_qS_l\end{eqnarray*}and identify a small (open)
disk  $Q\subset P_p$ about point $p$ with a subset of the Morse-chart $U_{x_k}$. (Such an identification is standard for a Euclidean space and
we have already identified $U_{x_k}$ with a ball in ${\bf R}^3$ via the Morse coordinate system $\eta_{x_k}$ (this identification also induces
scalar products on planes $P_p$ and $P_q$. As the double index $kl$ is fixed, we will supress it in the course of this proof.) Similarly, we
identify a disc $Q'\subset P_q$  about $q$ with a subset of  Morse chart $U_{y_l}$. We call the map: $$\zeta:Q\hookrightarrow  Q',\hspace
{.7cm}\zeta(v)=\lambda_v\cap Q'$${\it flow map}. By shrinking the disc $Q$ if necessary, we can suppose that the flow map is defined on $Q.
\vspace{.2cm}$

\noindent Note that the 2-plane $P_p$ is spanned by unit vectors $$E^\downarrow\in T_pS_k^\downarrow, \hspace{1cm}N^\downarrow\in( T_pS_k
^\downarrow)^\perp$$ choosen so that the pair $(T_p\lambda_p, E^\downarrow)$ is oriented in the sense of the fixed orientation on $W_k
^\downarrow$ while $(T_p\lambda_p, E^\downarrow,N^\downarrow)$ corresponds to the orientation of $M$. Similarly, we define the pair of vectors
$$(E\uparrow,N^\uparrow)\in T_qS_l$$.\vspace{.5cm}

{\noindent\bf  Definition 7.} The flow is {\it aligned along trajectory $\lambda$} if there exists an open neighborhood $Q^*\subset Q$
about point $p$ such that the flow map $$Q^*\to P_q$$ is a restriction of an orthogonal map $P_p\to P_q$ of 2-planes and its differential takes
vector $E^\downarrow$ into $N^\uparrow$ and it takes vector $N^\downarrow$ either into $E^\uparrow$ or into $-E^\uparrow$. We call such a set
$Q^*$ the {\it domain of alignment along trajectory $\lambda$}. Accordingly, the set $\zeta(Q^*)$ is the {\it range of alignment}. The flow is
{\it aligned} if  it is aligned along each trajectory with  index difference 1. $\vspace{.3cm}$

{\noindent}Let $\sigma: P_q\to P_q$ be the orientation preserving  linear transformation so that \begin{eqnarray*}\sigma\circ{\rm d}\zeta
(E^\downarrow)&=&N^\uparrow\\  \sigma\circ{\rm d}\zeta(N^\downarrow)&=&\pm E^\uparrow\end{eqnarray*}Observe that then $$\sigma\circ{\rm d}
\zeta:P_p\to P_q$$ is an orthogonal transformation.

{\noindent}Consider now maps $$\xi_1,\ \xi_2:Q'\hookrightarrow Q',\hspace{.5cm}\xi_1(v)=\zeta\circ{\rm d}\zeta^{-1}(v),\hspace{.5cm}\xi_2
={\rm id}_{Q'}.$$(we shrink $Q'$ when necessary). The theorem that is usually referred to as  "Equivalence of Tubular Neighborhoods" (see
e.g. \cite{Ma}) guarantees an isotopy $$H':Q'\times[0,1]\to Q'\hspace{1cm}{\rm with}\hspace{.3cm}{\rm supp}(H')\subset Q'$$ which is stationary at
point $q$ such that\begin{eqnarray}H'_1\circ\xi_1=\xi_2\circ\sigma\end{eqnarray}holds in a neighborhood of  $q$ (taken in $Q'$).

{\noindent}Lemma 4.7. in \cite{Mi} guarantees the existence of a gradientlike vector field $X'$ so that $$\zeta'=H'_1\circ\zeta$$ holds for
the new flow map $\zeta':Q\hookrightarrow Q'$. This means that$$\zeta'(v)=H'_1\circ\zeta(v)=H'_1\circ\xi_1\circ{\rm d}\zeta(v)
=\sigma\circ{\rm d}\zeta(v)\hspace{1cm}(v\in Q^*\subset Q)$$holds (in other words, {\underline{the new flow map coincides with an orthogonal
transforma-}}

{\noindent\underline{tion}})$$\sigma\circ{\rm d}\zeta|_{Q^*}:Q^*\to Q'$$ in a small neighborhood $Q^*\subset Q$, thus the flow is aligned
along trajectory $\lambda.\vspace{.3cm}$

{\noindent}From now on we presume that the flow is aligned. Let $Q_k^+$ denote the domain of alignment along trajectory $\mu_k^+$ (with
analoguos notation for alignment along $\mu_k^-$). We call the sets\begin{eqnarray*}{\cal A}^+&:=&\cup\{\lambda_p|\ p\in\mathop{\cup}\limits
_{k=1}^K Q_k^+\}\\ {\cal A}^-&:=&\cup\{\lambda_p|\ p\in\mathop{\cup}\limits_{k=1}^K Q_k^-\}\end{eqnarray*}the {\it regions of alignment}.
$\vspace{.3cm}$

{\noindent\bf Definition 8.} For a point $p\in S_k^\downarrow$ let $$v^+,\ v^-\in T_pS_k,\hspace{1cm}v^+,\ v^-\perp T_pS_k^\downarrow$$ be
the  two vectors such that:$\vspace{.2cm}$

$v^+\ (v^-)$ points at the component of $U_{x_k}\setminus(W_k^\downarrow\cup\{x_k\})$ that contains a

section of  trajectory $\mu^+\ (\mu^-)$ respectively.$\vspace{.3cm}$

Identifying the two vectors above with subsets of  $U_{x_k}$, they have maximal

length amongst vectors pointing to the same direction that are contained

in the region of alignment ${\cal A}^+\ ({\cal A}^-)$ respectively.$\vspace{.3cm}$

{\noindent}We call the two one parameter families of trajectories $${\cal M}_p^+:=\{\lambda_{tv^+}|\ t\in (0,1)\}\hspace{1cm}({\cal M}_p^-:=
\{\lambda_{tv^-}|\ t\in(0,1)\})$$ the {\it outbound, positively (negatively) oriented membranes} at point $p$. Trajectory $\lambda_p$ is
called the (common) {\it base} of these membranes. Note, that inclusions $${\cal M}_p^+\subset{\cal A}^+,\hspace{1cm}{\cal M}_p^-\subset{
\cal A}^-$$hold.\vspace{.5cm}

{\noindent\it  Observations:}$\vspace{.2cm}$

{\noindent\bf  I.:} To each membrane let's associate its base. Then we get a bijection between the set of positively oriented (outbound)
membranes and between set$${\cal W}^\downarrow:=\mathop{\cup}\limits_{k=1}^K{\cal W}_k^\downarrow$$A trajectory in ${\cal W}^\downarrow$
terminates at critical point $y_l$ iff each trajectory in the corresponding membranes terminate at $y_l$ (this follows from alignment). A
trajectory in ${\cal W}^\downarrow$ terminates at critical point $b_j$ iff each trajectory in the corresponding membranes terminate at $b_j.
\vspace{.2cm}$

{\noindent\bf II.:} Whenever all 2-discs $\{Q_k^+,\ Q_k^-\ 1\leq k\leq K\}$ have the same radius $\epsilon$, then (an outbound) membrane at a
point of  the outbound cirlcle $S_k^\downarrow$ is represented on the sphere $S$ by a geodesic section that initiates at one of
the points $x_k^+,\hspace{.2cm}x_k^-$ (depending upon the orientation of the membrane) and has length\begin{eqnarray}R:=3{\rm arctan}{
\epsilon\over\delta}\end{eqnarray} An application of  this observation to membranes consisting of trajectories that terminate at a
critical point of index one shows that for an aligned flow, an initial and a terminal  section of  a curve $\alpha_l^i$ is a
geodesic arc (section of a main circle on the sphere $S$). We parametrize the curves $\alpha_l^i$ so that these arcs are of form $$\alpha
_l^i((0,t_0))\hspace{.5cm}{\rm and} \hspace{.5cm}\alpha_l^i((1-t_0,1))$$ for each pair of indices $l=1,...,L,\ i=1,...,s_l.\vspace{.2cm}$

{\noindent\bf  III.:} We can define {\it inbound membranes} at points of  the circle $S_l^\uparrow$ in a similar fashion  (except that we do
not have the regions of alignment, that is, the sets that would correspond to sets  ${\cal A}^+$ and ${\cal A}^-$; so we just need to assign
some length ($<\delta'$) to the vectors $v^+$ and $v^-$ in the definition of a membrane. We will choose this length so, that whenever the
base point of vector $v^+\ (v^-)$ is close to some trajectory $\lambda_l^i$, then they are contained in the range of alignment along $\lambda
_l^i $. This choice ensures that - provided the base point is close enough to trajectory $\lambda_l^i $ - the inbound membranes are represented
on sphere $S$ by sections of circles that are centered at one of  the points $x_k^+\ (x_k^-)$ and have  (spherical) radius smaller than $R$.
The initial point of these sections belongs to one of  the arcs$$\alpha_l^i((0,t_0)),\hspace{1cm}\alpha_l^{i-1}((1-t_0,1)).\vspace{.3cm}$$

{\noindent}Let's parametrize the curve\begin{eqnarray}\alpha_l:(0,s_l)\hookrightarrow S\hspace{.7cm}\alpha_l(t)=\alpha_l^i(t-i+1)\hspace
{1cm}(i-1<t<i;\ i=1,...,s_l)\end{eqnarray}

Patching together the perpendicular geodesic tubular neighborhood of  curve $\alpha_l$ (taken in the
sphere $S$) with the family of sections on circles about points $x_k^\pm$ (that we described in {\bf III.}) yields a tubular neighborhood
\begin{eqnarray}(E_l,\rho_l,\alpha_l,O_l,\epsilon_l,\xi_l)\end{eqnarray} about each curve $\alpha_l\hspace{.2cm}
(1\leq l\leq L)$. Here the bundle $E_l\mathop{\longrightarrow}\limits^{\rho_l}\alpha_l$ is the perpendicular complement  of tangent bundle
$T\alpha_l\mathop{\longrightarrow}\limits^{\rho_l}\alpha_l$ in the bundle $T_{\alpha_l}S\mathop{\longrightarrow}\limits^{\rho_l}\alpha_l$.
The set $O_l$ is an open neighborhood of  $\alpha_l$ in $S,\hspace{.2cm}\epsilon_l:\overline\alpha_l\to[0,1)$ is a continuous function that
vanishes only at points in $\overline\alpha_l\setminus\alpha_l$ and $$\xi_l:E_l(\epsilon_l)\to O_l$$ is a diffeomorphism between  the
$\epsilon_l$-disc bundle of  $E_l$ and open set $O_l$. Observe that the fibers of  vector bundle $E_l$ are also parametrized  by $$E_l(t):=
E_l|_{\alpha_l(t)}\hspace{1.5cm}(t\in(0,s_l)\setminus{\bf N})$$

{\noindent\it Note:} We can suppose that$\vspace{.2cm}$

(i) The open neighborhoods $\{O_l|\ 1\leq l\leq L\}$ are pairwise disjoint.$\vspace{.3cm}$

(ii) The image $$O_l(t):=\xi_l(E_l(\epsilon_l)|_t)$$

is a connected open subset of a circle contained by sphere $S.\vspace{.3cm}$

(iii)  The two components of  the punctured sections $O_l(t)\setminus\alpha_l(t)$ are of

equal length so we can define\begin{eqnarray*}R_l:\alpha_l\to[0,1),\hspace{1cm}R_l(p)={\rm length\ of\ curve\  in\ fiber}\ (O_l\setminus
\alpha_l)|_p.\end{eqnarray*}

($R_l$ agrees with $\epsilon_l$ when restricted to the mid-part $\alpha_l^i((t_0,1-t_0))).\vspace{.3cm}$

It is not hard to see that - by using the "Equivalence of Tubular Neighborhoods" theorem - the flow can be re-arranged so that a
pair of inbound membranes are represented (on $S$) by the two components of the fiber $O_l(t)\setminus\alpha_l(t)$. As we saw above, this
characterization of inbound membranes already holds at points of  $\alpha_l^i(t)$ with $t\in(0,t')\cup(1-t',1)$ for some $t'<t_0$, so under
this rearrangement we can keep the flow unchanged in a neighborhood about trajectories with index difference 1.$\vspace{.6cm}$

The construction of the canonical characteristic map will be based on homotopy $H$ (see formula (3)) defined as follows: choose a
differentiable function$$g:(0,1)\to({1\over4},1)$$such that \begin{eqnarray*}g(\tau)&=&\tau+{1\over4}\hspace{1cm}0<\tau<{1\over 8}\\  g(\tau)
&=&\tau\hspace{1cm}{7\over8}<\tau<1\\ g(\tau)&>&\tau\hspace{1cm}0<\tau<{7\over 8}\\ {d\over d\tau}g&>&0\end{eqnarray*}For radius $R$ given in
formula (5) define homotopy\begin{eqnarray*}h_R&:&(0,R)\times[0,1]\to(0,R)\\  h_R(\tau,t)&=&tRg({\tau\over R})+(1-t)\tau\end{eqnarray*}

Let $D_k^+ \subset S$ be the geodesic disk centered at $x_k^+$ with (spherical) radius $R$. The disc $D_k^-\subset S$ about $x_k^-$ is defined
analoguosly (we suppose that $R$ is choosen so small, that all such discs are pairwise disjoint). Denote $\tilde x_k^+:={1\over3}x_k^+$ and
let ${\bf C}D_k^+$ be the (open) cone over  $D_k^+$ with vertex $\tilde x_k^+$ (i.e. $${\bf C}D_k^+:=\mathop{\cup}\limits_{x\in D_k^+}s_x$$
where $s_x:=\overline{x\tilde x_k^+}$ is the straight section that connects $x\in D_k^+$ with $\tilde x_k^+$).\vspace{.3cm}

{\noindent}First we define the restriction$$H|_{B_o\times[0,1]}.$$The restriction of homotopy $H$ to geodesic section $$\gamma:(0,R)\to D_k
^+\hspace{.5cm}({\rm with\ initial\  point}\  x_k^+)$$is defined as \begin{eqnarray}H|_\gamma(\gamma(\tau),t)=\gamma(h_R(\tau,t))\hspace{1cm}
(0<\tau<R,\ 0\leq t\leq1)\end{eqnarray} Define  $$H|_{{\bf C}D_k^+\times[0,1]}:{\bf C} D_k^+\times[0,1]\to{\bf C} D_k^+$$ by spherically
coning down the $H|_\gamma$'-s  (in  other words:\vspace{.2cm}

- $H$ is stationary at point $\tilde x_k^+.\vspace{.2cm}$

- $H_t(s_x)=s_{H|_\gamma(x,t)}\hspace{.5cm}(x\in\gamma,\ t\in[0,1]).\vspace{.2cm}$

- $\|H_t(z)\|=\|z\|\hspace{.2cm}(z\in s_x).\vspace{.3cm}$

Homotopy $H$ is defined similarly for the cones ${\bf C}D_k^-\hspace{.3cm}(k=1,...,K)$ and it is stationary on the rest of 3-ball $B.\vspace
{.3cm}$

{\noindent\it Notation:} Let $C_k^+\ (C_k^-)\hspace{.2cm}(1\leq k\leq K)$ denote  the circle on sphere $S$ with (geodesic) radius $R\over 4$,
centered at point $x_k^+\ (x_k^-)$ respectively and  let ${\bf  C}C_k^+\ ({\bf C}C_k^-)$  be the cone over the circle $C_k^+\ (C_k^-)$ with
vertex $\tilde x_k^+\ (\tilde x_k^-).\vspace{.3cm}$

{\noindent}Note that:$${\rm Bd}(H_1(B_o))\setminus S=\mathop{\cup}\limits_{k=1}^K({\bf C}C_k^+
\cup\{\tilde x_k^+\}\cup {\bf C}C_k^-\cup\{\tilde x_k^-\})\cup(\overline{H_1(B_o)}\setminus B_o)$$ Let
$$\tilde\alpha_l:={2\over3}\alpha_l$$ be the image of curve $\alpha_l$ under the dilation of ${\bf R}^3$ by factor $2\over3.\vspace{.4cm}$

{\noindent}For a fixed parameter $t\in(0,s_l)\setminus{\bf  N}$ consider the cone $${\bf  C}O_l(t):=\mathop{\cup}\limits_{x\in O_l(t)}
\overline{{\bf 0}x}$$over fiber $O_l(t)$ and with the origin $\bf 0$ as vertex. {\underline{Within each such cone}} take a "warped cone"$$
{\bf C}_l(t)\subset{\bf C}O_l(t)$$ over the base $O_l(t)\setminus\alpha_l(t)$ with vertex $\tilde\alpha_l(t)$. For parameter values$$i+t_0<t
<i+1-t_0\hspace{.3cm}(i=0,...,s_l-1)$$the cone ${\bf C}O_l(t)$ is contained in a 2-plane and the "warped cone" is just the regular one (with
vertex $\tilde\alpha_l(t)$). For a parameter value $t\in[i-t_0,i)\cup(i,i+t_0]$ and point $p\in O_l(t)$ let $c_p$ be the intersection
of  ${\bf C}O_l(t)$ with the plane through points $\tilde\alpha(t)$ and $p$ that is perpendicular to the fiber $O_l(t)$ at point $p$. Observe
that warped cones ${\bf C}_l(t)$ and ${\bf C}_l(t')$ are pairwise disjoint ($t\not=t'$). Let$${\bf C}_l:=\mathop{\cup}\limits_{t\in(0,s_l)
\setminus{\bf N}}{\bf C}_l(t).$$ We define homotopy $$H_l^*:{\bf  C}_l\times [0,1]\to{\bf C}_l$$ so that it keeps the fibers ${\bf C}_l(t)$
and in each such fiber the construction of  the homotopy will be analoguos to that of  $H|_{{\bf C}D_k^+\times[0,1]}$, namely:\vspace{.5cm}

(i) Parametrize the two components of  the set $$O_l(t)\setminus \alpha_l(t)$$

by archlength and, with $R=R_l(t)$,  use formula (8) to move each of  these

components within themselves.$\vspace{.3cm}$

 (ii) Cone down the above movement from point $\tilde\alpha_l(t)$, using curves $c_p$

$(p\in O_l(t))$ in place of sections $s_x.\vspace{.5cm}$

{\noindent}Define$$H^*:B_o\times[0,1]\to B_o$$as $H_l^*$ when restricted to the subset ${\bf C}_l\hspace{.2cm}(l=1,...,L)$ and stationary
elsewhere. Finally, let \begin{eqnarray}H_t=H^*_{t-1}\circ H_1\hspace{1cm}(1<t\leq 2)\end{eqnarray}

{\noindent}Let $$B^*:=H_2(B_o)\hspace{1.5cm}S^*:={\rm Bd}(H_2(B_o))$$
Then each ray in ${\bf R}^3$ that emanates from the origin intersects the boundary   $S^*$ in exactly one point, consequently $S^*$ is
homeomorphic to the 2-sphere.\vspace{.5cm}

{\noindent\it Observe that:}\vspace{.2cm}

$$H|_{\gamma t}(\alpha_l\cap D_k^+)\subset \alpha_l\cap D_k^+\hspace
{.3cm}(t\in[0,1],\ \alpha_l\cap D_k^+=\gamma)$$

{\noindent\it Notations:} Let\begin{eqnarray}q_l^{i*}&:=&\alpha_l^i((0,t_0))\cap\mathop{\cup}\limits_{k=1}^K(C_k^+\cup C_k^-)\\  q_l^{i**}&:=&\alpha
_l^i((1-t_0,1))\cap\mathop{\cup}\limits_{k=1}^K(C_k^+\cup C_k^-)\end{eqnarray}Note that there is a unique circle $C_k^+$ (or $C_k^-$) which
intersects the curve $\alpha_l^i((0,t_0))$ non-trivially. Let  $r_l^{i*}$ be the section between $q_l^{i*}$ and the vertex of the corresponding cone
${\bf C}C_k^+$ or ${\bf C}C_k^-$. Similar definition for $r_l^{i**}.\vspace{.4cm}$

For a fixed integer $1\leq l\leq L$ and value $t\in(0,s_l)\setminus{\bf  N}$ let $p^1_l(t),\ (p^2_l(t))$ be the points of  fiber $O_l(t)$ that
belong to parameter value $R_l(t)\over4$ (choosen so, that - taking into account the orientation of sphere  $S$ - point $p_l^1(t)$ resides on
the positive side of  curve $\alpha_l$). Let$${\bf  C}p_l^1(t):=c_{p_l^1(t)},\hspace{1cm}{\bf C}p_l^2(t):=c_{p_l^2(t)}$$and$${\bf  C}p_l^1
:=\mathop{\cup}\limits_{t\in(0,s_l)\setminus{\bf N}}{\bf  C}p_l^1(t)\hspace{1cm}{\bf  C}p_l^2:=\mathop{\cup}\limits_{t\in(0,s_l)\setminus{\bf
N}}{\bf  C}p_l^2(t)$$\begin{eqnarray*}\Sigma_1&:=&H_1^*(\mathop{\cup}\limits_{k=1}^K({\bf C}C_k^+\cup{\bf C}C_k^-)\cap B_o)\\ \Sigma_2&:=&S^*
\setminus(S\cup\Sigma_1\cup\mathop{\cup}\limits_{k=1}^K\{\tilde x_k^+,\ \tilde x_k^-\})\end{eqnarray*}

{\noindent\underline{Extension of the characteristic map to the boundary $S^*$:}}$\vspace{.3cm}$

{\noindent}We should define a map$$\psi:\overline{B^*}\to M$$such that $$\psi|_{B^*}=\phi_o\circ H_2^{-1}$$ So basically we need to define
$\psi$ on the boundary $S^*$. Observe, that\begin{eqnarray}\Sigma_2=\mathop{\cup}\limits_{l=1}^L({\bf C}p_l^1\cup{\bf  C}p_l^2\cup\tilde
\alpha_l)\cap\overline{B^*}\end{eqnarray}moreover, that  $S^*$ can be written as the disjoint union \begin{eqnarray}S^*=(S\cap S^*)\cup^*
(\mathop{\cup}\limits_{k=1}^K\{\tilde x_k^+,\ \tilde x_k^-\})\cup^*\Sigma_2\cup^*\Sigma_1\end{eqnarray}

Note that both maps  $\phi_o$ (see Definition 6.) and $H_2$ have (unique) continuous extensions $\overline\phi_o,\hspace{.2cm}\overline
H_2$, respectively, to the non-singular part of sphere $S$: $\overline\phi_o$ takes each connected component of  $S\setminus\overline
{(B\setminus B_o)}$ into a critical point of index 0, while $\overline H_2$ takes such a component into a component of  int($S\cap
S^*$). This implies, that the map $\overline\phi_o\circ \overline H_2^{-1}$  takes a connected component of  int($S\cap S^*$) into a critical
point of index 0. Define \begin{eqnarray}\psi({\rm component\ K\ of}\ S\cap S^*)&=&\overline\phi_o\circ\overline H_2^{-1}({\rm int}(K))\\
\psi(\tilde x_k^+)=\psi(\tilde x_k^-)&=&x_k\\ \psi(\tilde\alpha_l\cap\overline{B^*})&=&y_l\end{eqnarray}(When $\partial M\not=\emptyset$, then
the image $\psi(p)$ of a point $p\in S\cap S^*$ is defined as the terminal point of the trajectory of vector
field $X$ that corresponds to the ray through $p$.)\vspace{.3cm}

We recall that the components of fiber $O_l(t)\setminus\alpha_l(t)$ represent a pair of (inbound) membranes at the point of $S_l
^\uparrow$ that corresponds to $\alpha_l(t)$. Taking this into account we define\begin{eqnarray}\psi|_{{\bf C}p_l^1\cap S^*}&:&{\bf  C}p_l^1
\cap S^*\to\nu_l^+,\hspace{.7cm} \psi(p)=z\ \ni f(z)=3-\|p\|\\ \psi|_{{\bf C}p_l^2\cap S^*}&:&{\bf C}p_l^2\cap S^*\to\nu_l^-,\hspace{.7cm}
\psi(p)=z\ \ni f(z)=3-\|p\|\end{eqnarray}In other words, $\psi$ maps each point $p$ of the cone ${\bf C}p_l^1\cap S^*$ into the point of
curve $\nu_l^+$ which is at the level of  $3-\|p\|$ (note, that $0<3-\|p\|<1$ holds for  points of the cone ${\bf  C}p_l^1$), with analogous
definition for the other cone.\vspace{.3cm}

In order to define $\psi$ on $\Sigma_1$ note, that $$r_l^{i*}\cap\overline{ B^*}=r_l^{i*}\cap(B_o\cup\tilde\alpha_l^i)=\{z\in
r^{i*}_l |\ 1<\|z\|\leq2\}$$ and $r_l^{i*}$ is contained in the boundary of  the union of the $H_2$-image of rays through points of curve
$\alpha_l^i((0,t_0)).\ \phi_o$ maps each such ray into a trajectory that terminates at critical point $y_l$ (the $\phi_o$-images of these rays
comprise an (inbound) membrane with base $\lambda_l^i$)   so define \begin{eqnarray}\psi|_{ r^{i*}_l\cap B_o}:r^{i*}_l\cap B_o\to\lambda_l
^i\end{eqnarray}such that$$ f(\psi(p))=3-\|p\|$$ holds for $p\in r^{i*}_l\cap B_o$. We define map $\psi$ on $r^{i**}_l\cap B_o$ analogously
(it gets mapped onto trajectory $\lambda_l^{i+1})$.

{\noindent}Consider now a section $s_x$ for $x\in C_k^+\setminus\mathop{\cup}\limits_{l=1}^L\alpha_l$. The geodesic section $\gamma\subset
D_k^+$ through point $x$ represents a membrane at a point $q\in S_k^\downarrow$. Define \begin{eqnarray}\psi|_{H_1^*(s_x)}:H_1^*(s_x)\to
\lambda_q\hspace{1cm}f(\psi(p))=3-\|p\|\end{eqnarray}$\vspace{1cm}$

{\noindent\underline{Proof of continuity of map $\psi$}}$\vspace{.3cm}$

{\noindent}Note that by Definition 6. we have $$f(\psi(x))=f(\phi_o\circ H_2^{-1}(x))=3-\| H_2^{-1}(x)\|=3-\|x\|\hspace{.5cm}(x\in B^*)
$$while by formulas (14)-(20) we get $$f(\psi(x))=3-\|x\|\hspace{1cm}(x\in S^*)$$This implies that \begin{eqnarray}f(\psi(p_n))\longrightarrow
f(\psi(p))\end{eqnarray}for a convergent sequence $$p_n\longrightarrow p\hspace{.8cm}(p_n,\ p\in\overline{B^*})$$thus,
for continuity of  map $\psi$  it is enough to prove that for a convergent sequence $$p_n\longrightarrow p\hspace{1cm}(p_n\in\overline{B^*},\
p\in S^*)$$the corresponding sequence of  trajectories  converge, i.e.$$\lambda_{\psi(p_n)}\longrightarrow\lambda_{\psi(p)}$$in the following
sense of convergence:\vspace{.3cm}

{\noindent\bf  Definition 9.} We call a finite sequence of  joining trajectories a {\it cascade of trajectories}.  The number of trajectories
in the sequence is called the {\it  length} of the cascade.  The length of  a cascade can't  exceed the dimension of the manifold (in
our case, 3).\vspace{.3cm}

{\noindent\bf  Definition 10.} We say that a sequence of trajectories $\lambda_n$ {\it converges} to the cascade $\hat\lambda_1,...,\hat
\lambda_j\ (j=1,2,3)$  if the sequence of points $\lambda_n\cap f^{-1}(c)$ converges to $\hat\lambda_i\cap f^{-1}(c)$ whenever $\hat\lambda
_i\cap f^{-1}(c)\not=\emptyset\ (1\leq i\leq j).\vspace{.3cm}$

{\noindent}We prove continuity of map $\psi$ at points of  the boundary $S^*$ on a case-by-case basis. The cases correspond to the partition
of the sphere $S^*$ described in formula (13). $\vspace{.2cm}$

{\noindent} Consider a convergent sequence $$p_n\longrightarrow p\hspace{1.5cm}(p_n\in{\overline {B^*}},\ p\in S^*)\vspace{.3cm}$$

{\noindent\underline{Case 1. $p\in\Sigma_1$:}} Then $$(H_1^*)^{-1}(p)\in\mathop{\cup}\limits_{k=1}^K{\bf C}C_k^+\cup{\bf C}C_k^-$$So suppose
that $(H_1^*)^{-1}(p)\in{\bf C}C_k^+$, more precisely, that $(H_1^*)^{-1}(p)\in s_z$ for a point $z\in C_k^+$. Then $$(H_1^*)^{-1}(p_n)\in{\bf
C}D_k^+$$for large enough $n$, thus we have a sequence $$z_n\longrightarrow z\hspace{1.5cm}(z_n\in D_k^+)$$such that $$(H_1^*)^{-1}(p_n)\in
s_{z_n}$$Let $\gamma_n,\ \gamma\subset D_k^+$ be the geodesic sections (sections on main circles of  $S$ through point $x_k^+$) that contain
points $z_n,\ z$ respectively. The sequence of curves  $(\gamma_n)$ represent membranes at points of the outbound circle  $S_k^\downarrow$
and these membranes converge to the membrane represented by $\gamma$ (in the obvious sense). By formula (20), trajectory $\lambda_{\psi(p)}$
is the base of the membrane represented by $\gamma$ while trajectory $\lambda_{\psi(p_n)}$ is either contained in the membrane represented by
$\gamma_n$ or it is its base, consequently  $\lambda_{\psi(p_n)}\longrightarrow\lambda_{\psi(p)}.\vspace{.5cm}$

{\noindent\underline{Case 2. $p\in\{\tilde x_k^-,\ \tilde x_k^+|\ k=1,...,K\}$}}.  Suppose $p=\tilde x_k^+$. Then$$(H_1^*)^{-1}(p_n)
\longrightarrow\tilde x_k^+$$For the part of sequence $(p_n)$ that is contained in  $\Sigma_1$ the statement is an obvious consequence of
formula (20). For points $p_n\in B^*$ the sequence $H_2^{-1}(p_n)$ converges to $\tilde x_k^+$ thus the sequence of trajectories $\lambda
_{\psi(p_n)}$ should converge to $\mu_k^+\ \Rightarrow\hspace{.2cm}\psi(p_n)\longrightarrow\psi(p)=x_k.\vspace{.5cm}$

{\noindent}Case $p\in\Sigma_2$ splits into the following three subcases:\vspace{.3cm}

{\noindent\underline{Case 3. $p\in({\bf C}p_l^1\cup{\bf  C}p_l^2)\setminus\overline\Sigma_1$ for some $l=1,...,L$:}} we suppose $p\in{\bf C}
p_l^1\setminus\overline\Sigma_1$. For those members of sequence $(p_n)$ that belong to the boundary $S^*$ the image $\psi(p_n)$ is on
trajectory $\nu_l^+$, which is also the trajectroy through point $\psi(p)$ (see formula (17)). This means that it is enough to consider the
case when $p_n\not\in S^*$. Then $p_n\in{\bf C}_l$ for $n$ large and $p\in{\bf C}_l(t)$ for some parmeter value $t$. Observe, that the
sequence of rays (taken in ${\bf R}^3$) through points $(H_1^*)^{-1}(p_n)$ converges to the ray through $\alpha_l(t)$, thus the rays through
points \begin{eqnarray}H_1^{-1}\circ(H_1^*)^{-1}(p_n)\end{eqnarray} will converge to a ray through some $\alpha_l(t')\ (t'=t$ is possible)
moreover rays in formula (22) stay on the positive side of curve $\alpha_l$. This implies that the sequence of trajectories through
points$$\psi(p_n)=\phi_o\circ H_1^{-1}\circ(H_1^*)^{-1}(p_n)$$ should converge to a cascade with $\nu_l^+$ as its last member.$\vspace{.3cm}$

{\noindent\underline{Case 4. $p\in\tilde\alpha_l\setminus\overline\Sigma_1$:}} Then $\psi(p)=y_l$ by formula (16). By formulas (17) and (18),
$\psi$ maps the subsequence of $(p_n)$ that contained in cone ${\bf C}p_l^1\cup{\bf C}p_l^2$ either to trajectory $\nu^-_l$ or $\nu_l^+$ thus
-by formula (21)- the images $\psi(p_n)$ converge to critical point $y_l$. As homotopy $H_1^*$ is stationary at point $p$, this implies that
$$(H_1^*)^{-1}(p_{n_i})\longrightarrow p$$for the subsequence $(p_{n_i})$ of points that don't belong to the set  ${\bf C}p_l^1\cup{\bf C}
p_l^2$, thus$$H_1^{-1}\circ(H_1^*)^{-1}(p_{n_i})\longrightarrow H_1^{-1}(p)\in\tilde\alpha_l$$This shows, that critical point $y_l$ is a
limit point of the sequence of trajectories through points $$\psi(p_{n_i})=\phi_o\circ  H_1^{-1}\circ(H_1^*)^{-1}(p_{n_i}).$$

{\noindent\underline{Case 5. $p\in\overline\Sigma_1\cap\Sigma_2$}}: Observe that each connected component of $\overline\Sigma_1\cap\Sigma_2$
is an arch (homeomorphic image of the $(0,1)$ interval). We call these connected components {\it edges}. Each curve $\tilde\alpha^i_l$ intersects
exactly two edges (non-trivially) that are denoted by $$e^{i*}_l,\hspace{.3cm} e^{i**}_l$$ respectively.  Each edge is contained in the boundary
of a  (uniquely determined) element of the set $$\{H_1^*({\bf C}C_k^+\cap B_o),\ H_1^*({\bf C}C_k^-\cap B_o)| 1\leq k\leq K\}$$

{\noindent}To prove continuity at a point $p\in e^{i*}_l$ we combine the two reasonings that were given in Cases 1. and 3. We decompose the
sequence $p_n\longrightarrow p$ into two subsequences\vspace{.2cm}

(i)  $p_{n_j}'\in\Sigma_1\cup  H_2(B_o) \vspace{.2cm}$

(ii)  $p"_{n_j}\in\Sigma_2$

{\noindent}For clause (i) we consider the sequence of geodesic sections $\gamma_{n_j}$ just as we did in Case 1. These geodesic sections
converge to the geodesic section $\alpha_l^i(0,t_0)$ and stay on one side of  $\alpha_l^i(0,t_0)$ consequently $\lambda_{\psi(p_{n_j}')}$
converges either to cascade $(\mu^+,\lambda_l^i,\nu^+)$ or to $(\lambda_l^i,\nu^+)$ and this was to be proven.

The discussion we gave in Case 3. applies to the situation in clause (ii) word-by-word.\vspace{.3cm}

{\noindent\underline{Case 6.  $p\in S$:}} as continuity of  map $\psi$ has already been established at interior points of the connected
components of  $S\cap S^*$, we can suppose that $p$ is a boundary-point of such a component $K$. When $\partial M\not=\emptyset$ then we prove
continuity at points of component $K$ by the technique of converging membranes (see cases 1-5). Otherwise $\psi({\rm int}(K))$ is a critical
point of index $0$ and Bd$(K)\subset(\overline\Sigma_1\cup\overline\Sigma_2)\cap S$. For a point $p\in$Bd$(K)$ case $p\in\overline\Sigma_1
\cap\overline\Sigma_2\cap S$ (i.e. when $p$ is the boundary point of some edge) is the most intricate one (and the other two cases$$p\in
(\overline\Sigma_1\setminus\overline\Sigma_2)\cap S\hspace{1cm}p\in(\overline\Sigma_2\setminus\overline\Sigma_1)\cap S$$can be treated similarly),
we presume that $p\in$Bd$(K)$ also belongs to the boundary of some edge. Then the image $\psi(p)$ is a critical point of index 0 and it is
not hard to check that $\psi(p)=\psi(K)$. Note that - by construction - point $p$ has an open neighborhood $U\subset\overline{B^*}$ such that
$U\cap S\cap S^*\subset K$. Using reductio ad absurdum consider a sequence $$p_n\longrightarrow p,\hspace{.2cm}(p_n\in U)\hspace{1cm}\psi(p_n)
\longrightarrow x\not=\psi(p)$$ Then, by formula (21), $x$ must be a critical point of index 0, thus the sequence $(p_n)$ can intersect
int($K$) in finitely many points only, thus we can suppose that
it does not intersect it at all. Observe, that by the choice of neighborhood $U:\vspace{.2cm}$

(i)  Trajectory $\lambda_{\psi(q)}$ terminates at critical point $\psi(p)$ for $q\in U\cap B^*.
\vspace{.2cm}$

(ii)  Each point $q\in U\cap(\Sigma_1\cup\Sigma_2)$ is a limit point of the set $U\cap B^*$ and

continuity of $\psi$ at $q$ has already been proven, thus we infer, that  $\lambda_{\psi(q)}$

should terminate at critical point $\psi(p)$ while $\lambda_{\psi(p_n))}$ is supposed to ter-

minate at $x$, contradiction.

{\noindent}Proof of continuity of map $\psi$ is now complete.\hfill\rule{3mm}{3mm}
 $\vspace{.5cm}$

{\noindent}Let $r_z$ denote the ray through $z\in{\bf R}^3\setminus{\bf 0}$ and define$$\rho:\overline{B^*}\to[1,3],\hspace{.8cm}\rho(z)=
\|r_z\cap S^*\|$$Then map\begin{eqnarray}z\longrightarrow{z\over\rho(z)}\end{eqnarray}is a homeomorphism between $\overline{B^*}$ and the
closed unit 3-ball. The inverse of this map composed with $\psi$ is only piecewise differentiable, but either by smoothing or by applying the more intricate technique of pathcing
together local coordinates (see e.g. in Major \cite{Maj}) we can make it smooth in the interior of the unit ball. This proves Theorem B.
\vspace{.5cm}

Only cluase (iii) of Theorem A. needs further consideration. We define the stratified conic structure about trajectory spaces as follows:
\vspace{.3cm}

(i) By using the Morse domains, let the fibers of the tubular neighborhood

about a critical point $z$ be straight sections of length $\delta$ that emanate from $z.\vspace{.2cm}$

(ii) A connected component of a 1-dimensional stratum is a trajectory $\lambda$ of

vector field $X$ (with index difference 1). Let the tube $O_\lambda$ about $\lambda$ be the part

of the region of alignment about $\lambda$ that falls between the levels of the end-points

of $\lambda$ (e.g. for a trajectory $\lambda\in{\cal W}_{kl}$ let the tube be$$f^{-1}(1,2)\cap\bigcup\{\lambda_q\ |\ q\in Q^*\}$$

where $Q^*$ is the domain of alignment about $\lambda$. For a vector $v\in T_pQ^*$ at point

$p\in\lambda\cap Q^*$ consider the 1-parameter family of trajectories$$\lambda_v(t):=\lambda_{tv}\ (t\in[0,1))$$

and at a point $q=\lambda\cap f^{-1}(c)$ let the fiber of the tubular neighborhood be

curve $\lambda_v(t)\cap f^{-1}(c)$ with tangent vector ${\mathop{\lambda}\limits^.}_v(0)$. Taking properties of the flow-

map into account it is not hard to see, that this way we get a thick tubular

neighborhood about trajectory $\lambda$.\vspace{.2cm}

(iii) Similarly, about 2-dimensional stratum$$\mathop{\cup}\limits_{k=1}^KW_k^\downarrow\setminus\big(\mathop{\cup}\limits_{l=1}^LW_l^\uparrow
\big)$$

let the tube be$$(f<2)\cap({\cal A}^+\cup{\cal A}^-\cup\mathop{\cup}\limits_{k=1}^KW_k^\downarrow)\setminus\mathop{\cup}\limits_{l=1}^LW_l^\uparrow$$

and define the fibers of the tubular neighbrohood as the intersection of

mebranes with level sets of the Morse function.\vspace{.2cm}

(iv) The tubular neighborhood of a trajectory space connecting $o$ with a

critical point of index $1$ is defined in complete analogy with (iii).\vspace{.5cm}

{\noindent}It is not hard to see that formula (1) holds for the retractions, proving Theorem A. this way.\vspace{.5cm}

As for the generalization of our theorems to higher dimensions let's note, that the fundamental difficulty is to produce (or even: define) the aligned
flow. Once the flow is aligned along all trajectory spaces of vector field $X$, the successive isotopies (i.e. the push-in's) of the $n$-ball
can be defined in complete analogy with our isotopy $H$ (see formula (3)) and then the characteristic map $\psi$ can also be defined in a
similar way; furthermore the proof of its continuity would not require any novel steps either. We plan to discuss alignment in higher dimensions
in an upcoming paper. The major goal here is to prove the following generalization of clause (ii) of Theorem A. (stated below as a conjecture)
\vspace{.4cm}

{\noindent\bf Conjecture} Let $(f,X)$ be a Morse-pair on a compact orientable n-dimensional manifold $M$ where $f$ has only one critical point
$o$ of index $n$. Then there is a diffeomorphism of $M$ (isotopic to id$_M$ and the identity on a neighborhood of the critical set) which
transforms $(f,X)$ so that for the trajectory spaces of the new Morse-pair $(f',X')$ the following holds:\vspace{.2cm}

 (i) There is a continuous surjection $\psi$ from the closed unit ball onto $M$ such

 that its restriction to the open ball is a diffeomorphism onto $n$-cell $W_o^\downarrow.$\vspace{.2cm}

(ii) For each critical point $x\in$Cr$(f)$ of index $0<i<n$ the pre-image

$A_x=\psi^{-1}(x)\subset S^{n-1}$ is a submanifold of dimension $n-i-1$ and$$\psi({\bf C}A_x)=W_{ox}$$]

holds (where ${\bf C}A_x$ is the cone over $A_x$ with vertex $\bf 0$).\vspace{.2cm}

 (iii) About each such submanifold $A_x\subset S^{n-1}$ there exists a tubular neighbor-

 hood so that $\psi$ maps each fiber of bundle $\xi_x(E_x(\epsilon_x))$ diffeomorphically onto

 the submanifold $W_x^\downarrow\cup\{x\}\subset M$.\vspace{.5cm}

In his paper (see \cite{V1}, Proposition 2.6.) Verona considers a prestratification of an $n$-dimensional manifold $A$ and defines a closed
"tubular neighborhood" about the closure of the codimension-1 stratum, with the boundary $W$ of the tube as base, so that $W\subset A$ is a
smooth $n-1$-dimensional submanifold. His conept fails to be a tubular neighborhood in the usual sense, since his "$\xi$" (he denotes it othwerwise)
maps $W\times[0,1]\to A$ so that:\vspace{.3cm}

 $\xi|_{W\times(0,1]}:W\times(0,1]\to A$ is a smooth imbedding

 $\xi(W\times\{0\})$ is the closure of the codimenion-1 stratum\vspace{.5cm}

{\noindent}(so the restriction  $\xi|_{W\times\{0\}}:W\times\{0\}\to A$ is not necessarily  injective.)\vspace{.4cm}

In our special case it is not hard to construct such a surface $W$: the images of the fibers we defined in the proof of clause (iii) of
Theorem A. under the map $H_2\circ\phi_o^{-1}$ yield a tubular neighborhood (in closed ball $\overline{B^*}$) of the differentiable part of
2-sphere $S^*$ (supplemented by small sections on rays over part $S\cap S^*$). Further taking the image of these fibers by the map given in
(23) is a part of a tubular neighborhood of the unit 2-sphere $S$ (taken in the closed unit 3-ball), with some curves in $S$ over which the
bundle is not yet defined. It is not hard to supplement the bundle overe these curves, and then smoothen it. This way we get a tubular
neighborhood of $S^2\subset\overline{B^3}$. The $\psi$-image of a section of this tubular neighborhood is an appropriate 2-surface (sphere) $W$
described in Verona \cite{V1}.

In the light of the above said, it is not hard to see, that conic stratification and Verona's $W$ grab two sides of the same thing:\vspace{.3cm}

- Either we take a "global" submanifold $W$ and then, for $W\times[0,1]$ we collapse

certain submanifolds of $W\times\{0\}$\vspace{.2cm}

- We take a tubular neighborhood about each stratum (i.e. trajectory space of

vector field $X$) and then stipulate, that they fit together nicely (i.e. that

they obey the rules of a conic stratification).\vspace{.4cm}

Finally let's remark, that the trajectory spaces of an aligned flow join smoothly (i.e. at the degree of differentiability of the flow) instead
of being only $C^1$ (the way found in Laudenbach [1]).

\vspace{1cm}

{\noindent}IMRE MAJOR, CENTRAL EUROPEAN UNIVERSITY and DENES GABOR COLLEGE OF INFORMATICS, BUDAPEST\vspace{.5cm}

e-mail: imajor@hotmail.com


\begin{thebibliography}{99}

\bibitem{L}Laudenbach, F.: {\it On the Thom-Smale complex.}

(Appendix to Bismut, J-M., Zhang, W.: An extension of a theorem by Cheeger and Muller. Asterisque {\bf 205} (1992)

\bibitem{Maj}Major, I.: Can one tell if a presentation describes the fundamental group of an irreducible 3-manifold? (preprint).

\bibitem{Ma}Mather, J.: {\it Notes on topological stability}. Harvard University, 1970.

\bibitem{Mi1}Milnor, J.W.:  {\it Morse Theory}.

Ann. of Math. Studies, Princeton Univ. Press, (1963).

\bibitem{Mi}Milnor, J.W.: {\it Lectures on the H-cobordism theorem}.

Ann. of  Math. Studies {\bf 51}, Princeton Univ. Press, (1965).

\bibitem{V1} Verona, A.: {\it Homological properties of abstract stratifications.}

Rev. Roumanie Math. Pures Appl. Vol. XVII. No. 7. pp. 1109-1121 (Bucharest, 1972)

\bibitem{V2} Verona, A.: {\it Stratified mappings - structure and trianguability.}

LNM, Springer Verlag, 1984.


\end{thebibliography}
\end{document}